\documentclass[11pt]{gtart}
\usepackage{amsmath, amssymb, graphicx, epsfig, verbatim}
\usepackage{epic}

\newtheorem{thm}{Theorem}[section]

\newtheorem{cor}[thm]{Corollary}
\newtheorem{lem}[thm]{Lemma}
\newtheorem{prop}[thm]{Proposition}

\theoremstyle{definition}

\newtheorem{rem}[thm]{Remark}

\numberwithin{equation}{section}

\newcommand{\bfz}{{\mathbb {Z}}}

\newcommand{\Z}{\mathbb Z}

\newcommand{\cp}{{\mathbb C}{\mathbb P}}
\newcommand{\cpkk}{{\overline {{\mathbb C}{\mathbb P}^2}}}
\newcommand{\cpk}{{\mathbb {CP}}^2}
\newcommand{\cphat}{{\mathbb {CP}}^2\# 6{\overline {{\mathbb C}{\mathbb P}^2}}}
\newcommand{\cpot}{{\mathbb {CP}}^2\# 5{\overline {{\mathbb C}{\mathbb P}^2}}}
\newcommand{\cphet}{{\mathbb {CP}}^2\# 7{\overline {{\mathbb C}{\mathbb P}^2}}}
\newcommand{\cpnyolc}{{\mathbb {CP}}^2\# 8{\overline {{\mathbb C}{\mathbb P}^2}}}
\newcommand{\eegy}{{\mathbb {CP}}^2\# 9{\overline {{\mathbb C}{\mathbb P}^2}}}

\begin{document}

\title{Exotic smooth structures on $\cpot$}

\author{Jongil Park}
\address{Department of Mathematical Sciences\\
Seoul National University\\
San 56-1, Shinlim-dong, Gwanak-gu\\
Seoul 151-747, Korea}

\email{jipark@math.snu.ac.kr}

\author{Andr\'{a}s I. Stipsicz}
\address{R\'enyi Institute of Mathematics\\
Hungarian Academy of Sciences\\
H-1053 Budapest\\ 
Re\'altanoda utca 13--15, Hungary and\\
Institute for Advanced Study, Princeton, NJ}

\email{stipsicz@renyi.hu}

\secondauthor{Zolt\'an Szab\'o}
\secondaddress{Department of Mathematics\\
Princeton University,\\
 Princeton, NJ, 08544}
\secondemail{szabo@math.princeton.edu}

\date{Preliminary version}
\begin{abstract}
Motivated by a construction of Fintushel and Stern, we show that the
topological 4--manifold $\cpot$ supports infinitely many distinct
smooth structures.
\end{abstract}
\maketitle

\section{Introduction}
It is a basic problem in 4--dimensional topology to find exotic smooth
structures on rational surfaces. The first such structures were found by
Donaldson~\cite{D1}; these examples were homeomorphic to
$\eegy$. While in this homeomorphism type many exotic examples were
constructed \cite{FSknot, FM, Sz}, the cases of $\cpk \# k\cpkk$ with
$k<9$ were more elusive. The Barlow surface \cite{Bar} provided the
first exotic structure on $\cpk \# 8 \cpkk$, see \cite{Kot}.  More
recently, an exotic smooth structure on $\cphet$ has been constructed
\cite{P}.  After this example many new exotic 4--manifolds with small
Euler characteristic have been found.  In \cite{SSz} symplectic
4--manifolds homeomorphic but not diffeomorphic to $\cphat$ were
constructed, implying the existence of an exotic smooth structure on
$\cphat$. In a beautiful recent paper \cite{FSuj} Fintushel and Stern
showed the existence of infinitely many distinct smooth structures on
$\cpk \# k \cpkk$ with $k=6,7,8$. Combining their technique of knot
surgery in a double node neighborhood with a particular form of
generalized rational blow--down, in this note we prove

\begin{thm}\label{t:main}
There exist infinitely many pairwise nondiffeomorphic 4--manifolds 
all homeomorphic to $\cpot$.
\end{thm}

In Section~\ref{s:mas} various constructions of 4--manifolds
homeomorphic to $\cpot$ are described. In Section~\ref{s:harm}
we use Seiberg--Witten theory to show that many of these examples are
mutually nondiffeomorphic, leading us to the proof of 
Theorem~\ref{t:main}.

{\bf Acknowledgement:} We would like to thank Ron Fintushel and Ron
Stern for sending us an early version of \cite{FSuj} which provided an
essential ingredient to the construction given in this paper.
We also thank Mustafa Korkmaz for useful e--mail exchange. JP
was supported by the Korea Research Foundation Grant (KRF-2004-013-C00002), 
AS by OTKA T49449 and ZSz by NSF grant number DMS 0406155.

\section{The constructions}
\label{s:mas}
We construct our examples using knot surgery (in a double node
neighborhood, as in \cite{FSuj}) when applied to particular elliptic
fibrations. The special properties of the chosen elliptic fibration
allow us to find a configuration in the result of the knot surgery
such that after rationally blowing it down we arrive to a 4--manifold
homeomorphic, but not diffeomorphic to $\cpot$. By using a suitable
infinite set of knots (the twist knots already encountered in
\cite{FSuj}, cf. also \cite{FSknot, Sz}), we get an infinite family of
4--manifolds all homeomorphic to $\cpot$.

\subsection{Elliptic fibrations}
Singular fibers of holomorphic elliptic fibrations have been
classified \cite{Kod} (cf. also \cite{HKK}). In this note we will
consider fibrations containing only singular fibers of type $I_n$
($n\geq 1$). Recall that the singular fiber $I_1$ (also known as the
\emph{fishtail} fiber) is an immersed 2--sphere with one positive double
point, and it is created from a regular torus fiber by collapsing a
homologically essential simple closed curve (the {\em vanishing cycle}
of the singular fiber).  The $I_n$--fiber ($n\geq 2$) is a collection
of $n$ 2--spheres of self--intersection $(-2)$, intersecting each
other in a circular pattern, see \cite[page~35]{HKK}. An elliptic
fibration with singular fibers only of type $I_n$ are {\em Lefschetz
fibrations} in the sense of \cite[Chapter~8]{GS}. The only subtlety we
have to keep in mind is that here we allow a singular fiber to contain
more than one singular points as well.

Lefschetz fibrations can be conveniently described by the monodromy
factorization induced by the singular fibers of the fibration, that is, by 
a word involving right--handed Dehn twists
which is equal to 1 in the mapping class group of the regular fiber.
The mapping class group $\Gamma _1$ of the 2--torus $T^2$ can be presented
as 
\[
\Gamma _1 =\{ a,b\mid aba=bab, (ab)^6=1\} ,
\]
where $a,b\in \Gamma _1$ denote the right--handed Dehn twists along
the two standard simple closed curves $A,B$ in $T^2$ intersecting each
other transversally in a unique point.
This group can identified with SL$(2;\Z)$ by mapping $a$ to 
$\left(\begin{smallmatrix}
 1 & 1 \\
 0 & 1 
\end{smallmatrix}\right)$ and $b$ to
$\left(\begin{smallmatrix}
 1 & 0 \\
 -1 & 1 
\end{smallmatrix}\right)$.
For example, the standard elliptic fibration we get by blowing up nine
base points of a generic elliptic pencil in $\cpk$ results the
monodromy factorization $(ab)^6$. Using the braid relation $aba=bab$
it can be shown that $(a^3b)^3$ also defines an elliptic fibration on
$\eegy$.  Furthermore, it is easy to see that for any expression $x\in
\Gamma _1$ the mapping class $a^x=xax^{-1}$  can be
identified with the right--handed Dehn twist along the image of $a$
under a map giving $x$. Note, for example, that the braid relation 
implies that $b=a^{ab}$.

The monodromy of a fishtail fiber can be shown to be equal to the
right--handed Dehn twist along the vanishing cycle corresponding to
the given singular fiber. An $I_n$--fiber can be created by collapsing
$n$ parallel (homologically essential) simple closed curves, therefore
the monodromy of such a fiber is equal to the $n^{\rm {th}}$ power of
the right--handed Dehn twist along one of the parallel curves.

In our constructions we will need the existence of a section, which can 
also be read off from the monodromy factorization. In general, a Lefschetz
fibration admits a section if the monodromy factorization induced by it 
can be lifted from the mapping class group of its generic fiber to the
mapping class group of the fiber with one marked point. In the case 
of a genus--1 Lefschetz fibration, however, the forgetful map $f\colon
\Gamma _1 ^1 \to \Gamma _1$ mapping from the mapping class group $\Gamma _1^1$
of $T^2$ with one marked point to $\Gamma_1$ is an isomorphism, implying
in particular
\begin{lem}\label{l:szeles}
Any genus--1 Lefschetz fibration over $S^2$ admits a section. \qed
\end{lem}

\subsection{The definition of the 4--manifold $X_n$}

\begin{figure}[htb]
\begin{center}
\setlength{\unitlength}{1mm}
\includegraphics[height=7cm]{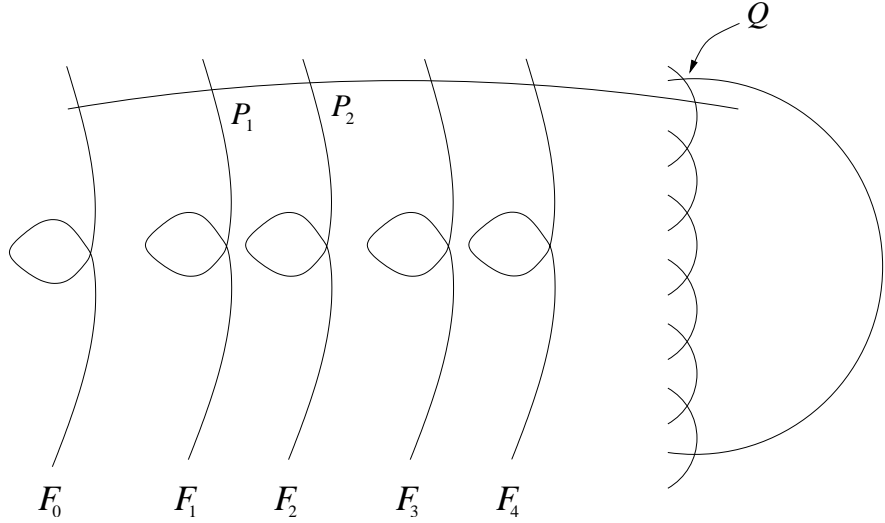}
\end{center}
\caption{\quad The schematic diagram of the fibration with an $I_7$ fiber}
\label{f:I7}
\end{figure}

Our first construction of exotic 4--manifolds relies on the following
existence result.  (For a schematic picture of the fibration see
Figure~\ref{f:I7}.)

\begin{prop}\label{p:letI7}
There exists an elliptic fibration $\eegy \to \cp ^1$ with five
fishtail fibers, an $I_7$--fiber and a section. Furthermore, we can
assume that two of the five fishtail fibers have isotopic vanishing
cycles. 
\end{prop}
\begin{proof}
We will show the existence of such fibration by finding an appropriate
factorization of $1$ in the mapping class group $\Gamma _1$ of the
torus.  Start with the fibration on $\eegy$ defined by the
factorization
\[
(a^3b)^3
\]
of $1\in \Gamma _1$. Notice that
\[
(a^3b)^3=a^7(a^{-1}(a^{-3}ba^3)a)(a^{-1}ba)a^2b=
a^7b^{a^{-4}}b^{a^{-1}}a^2b .
\]
Since $a^7$ is the monodromy of an $I_7$--fiber, its existence in the
above fibration is verified. The term $a^2$ gives rise to two fishtail
fibers with isotopic vanishing cycles in the complement of the
$I_7$--fiber.  Finally, Lemma~\ref{l:szeles} shows the existence of a
section in the fibration.
\end{proof}

Suppose now that $p>q>0$ are relatively prime integers.
Let us define the 4--manifold $C_{p,q}$ as the result of the linear
plumbing with weights specified by the continued fraction coefficients
of $-\frac{p^2}{pq-1}$. It is known \cite{CH} (cf. also
\cite{FS1, Pratb, SSz, Sym}) that the boundary 
$\partial C_{p,q}=L(p^2,pq-1)$ is a lens space which bounds a rational ball
$B_{p,q}$. The replacement of  an embedded copy of $C_{p,q}\subset X$ with 
$B_{p,q}$ is called the {\em (generalized) rational blow--down} of
$X$ along $C_{p,q}$. This operation was introduced and successfully 
applied by Fintushel and Stern \cite{FS1} in the case of $q=1$ and 
studied in \cite{Pratb, Sym} in the above generality.

Now we are ready to turn to the construction of the 4--manifolds
homeomorphic but not diffeomorphic to $\cpot$.  Let $K_n$ denote the
$n$--twist knot as it is depicted in \cite{FSuj}.  Let $F_3, F_4$ of
Figure~\ref{f:I7} denote the fishtail fibers with isotopic vanishing
cycles.  Following the convention of \cite{FSuj} we denote the result
of the knot surgery in a double node neighborhood containing $F_3,F_4$
and with knot $K_n$ by $Y_n$.  Fintushel and Stern \cite{FSuj} prove
the existence of a ``pseudo--section'' $S\subset Y_n$ which is an
immersed sphere with one positive double point, homological self--intersection
$-1$, and which transversally intersects $F_1,F_2$ and one of the
spheres in the $I_7$--fiber: The section of the fibration, punctured
by the fiber along which the knot surgery is performed, can be glued
to the genus--1 Seifert surface of the knot $K_n$.  In this way an
embedded torus $T$ of self--intersection $-1$ is found in $Y_n$.
Using the two thimbles of the isotopic vanishing cycles, Fintushel and
Stern find a disk $U$ attached to $T$ with relative self--intersection
$-1$.  {From} $T$ and $U$ now it is an easy task to find the immersed
sphere with a positive double point and which is homologous to
$T$. For more details of the construction see \cite{FSuj}.

\begin{figure}[htb]
\begin{center}
\setlength{\unitlength}{1mm}
\includegraphics[height=6cm]{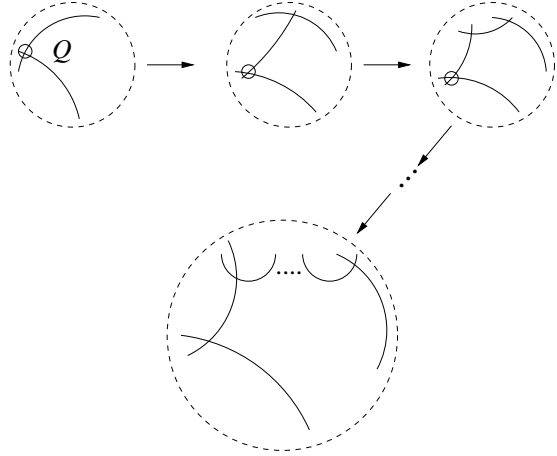}
\end{center}
\caption{\quad Infinitely close blow--ups at $Q$}
\label{f:infinite}
\end{figure}

Let us blow up
$Y_n$ in the double point of the pseudo--section, and in the double
points of the fishtail fibers $F_1$ and $F_2$.  After smoothing the
intersections $P_1,P_2$, we get a sphere of self--intersection $-9$
intersecting the $I_7$--fiber transversally at one point.  Now we apply
eight infinitely close blow--ups at the point $Q$ as it is  shown by
Figure~\ref{f:infinite}. 
This construction results in a chain of 2--spheres, with a neighborhood
diffeomorphic to the 4--manifold $C$ we get by plumbing along a linear
chain with weights
\[
(-9,-10,-2,-2,-2,-2,-2,-3,-2,-2,-2,-2,-2,-2,-2)
\]
in the eleven--fold blow--up of $Y_n$. 
Simple computation identifies $C$ with $C_{71,8}$.
Define $X_n$ as the (generalized) rational blow--down of $Y_n$
along $C$, that is, 
\[
X_n = (Y_n \# 11 \cpkk -{\mbox { int }}C) \cup B_{71,8} .
\] 

\begin{thm}
$X_n$ is homeomorphic to $\cpot$.
\end{thm}
\begin{proof}
The 4--manifold $Y_n$ has trivial fundamental group, since the
fibration admits a section and two different vanishing cycles in the
complement of the double node neighborhood.  Simple connectivity of
$X_n$ follows from the fact that the complement of $C$ in $Y_n \# 11
\cpkk$ is simply connected, since the generator of $\pi _1 (\partial
C) $ can be contracted along the fishtail fiber $F_0$ present in the
fibration but not used in constructing the configuration $C$ and from
the surjectivity of the natural map $\pi _1 (\partial B )\to \pi
_1(B)$.  Now simple Euler characteristic and signature computation
together with Freedman's Theorem on the classification of topological
4--manifolds \cite{Fr} imply the result.
\end{proof}

\subsection{Further constructions}
Many similar constructions can be carried out using different elliptic
fibrations or different sets of knots. Below we outline constructions
relying on various types of elliptic fibrations.

\subsubsection{Another construction using the $I_7$--fiber}
A similar argument provides an embedding of $C_{212,55}$ into $Y_n\#
12\cpkk$ by smoothing only at $P_2$ and keeping the transverse
intersection $P_1$. In this case one further blow--up of a
$(-2)$--sphere is necessary, leading to the chain
\[
(-4,-7,-10,-2,-2,-2,-2,-2,-3,-2,-2,-2,-2,-3,-2,-2)
\] 
in $Y_n\# 12\cpkk$.  Blowing this configuration down we get a sequence
of 4--manifolds with the same properties as $X_n$. (The hemisphere
originated from the exceptional sphere of the last blow--up can be
used to show that the resulting configuration of spheres in the
twelve--fold blow--up of $Y_n$ has simply connected complement.)

\subsubsection{Configurations using the $I_8$--fiber}
Many other examples can be given using the $I_8$--fiber. To see the
existence of the required fibration, we need a result similar to
Proposition~\ref{p:letI7}.
\begin{prop}
There exists an elliptic fibration $\eegy \to \cp ^1$ with four
fishtail fibers, an $I_8$--fiber and a section. Furthermore, we can
assume that two of the four fishtail fibers have isotopic vanishing
cycles. 
\end{prop}
\begin{proof}
Using the braid relation it is fairly easy to see that the expression
\[
a^3ba^2b^2a^2ba
\]
is equal to 1 in $\Gamma _1^1$, hence defines an elliptic fibration
with a section. 
Since it can be written as
\[
a^8(b^{a^{-2}})b^2(b^{a^2}),
\]
the resulting fibration can be chosen to have an $I_8$--fiber and two
fishtails in its complement with isotopic vanishing cycles.
\end{proof}
Our further constructions rely on
\begin{prop}\label{p:ujconf}
Let $Y_n$ be the 4--manifold defined above.
\begin{enumerate}
\item The 4--manifold $C_{44,9}$ embeds into $Y_n \# 8 \cpkk$;
\item $C_{79,44}$ admits an embedding into $Y_n \# 9\cpkk$;
\item $C_{89,9}$ embeds into $Y_n \# 13\cpkk$;
\item $C_{169,89}$ can be embedded into $Y_n\# 14 \cpkk$;
\item $C_{301,62}$ admits an embedding into $Y_n \# 14\cpkk$; and finally
\item $C_{540,301}$ is a submanifold of $Y_n \# 15 \cpkk$.
\end{enumerate}
The complements of these configurations are simply connected.
\end{prop}
\begin{rem}
Recall that the above 4--manifolds can be given by the linear
plumbings as follows:

\noindent $C_{44,9}=(-5,-11,-2,-2,-2,-2,-2,-2,-3,-2,-2,-2)$,

\noindent $C_{79,44}=(-2,-5,-11,-2,-2,-2,-2,-2,-2,-3,-2,-2,-3)$,

\noindent $C_{89,9}=(-10,-11,-2,-2,-2,-2,-2,-2,-3,-2,-2,-2,-2,-2,-2,-2,-2)$, 

\noindent
$C_{169,89}=(-2,-10,-11,-2,-2,-2,-2,-2,-2,-3,-2,-2,-2,-2,-2,-2,-2,-3)$,

\noindent
$C_{301,62}=(-5,-7,-11,-2,-2,-2,-2,-2,-2,-3,-2,-2,-2,-2,-3,-2,-2,-2)$
and finally

\noindent
$C_{540,301}=(-2,-5,-7,-11,-2,-2,-2,-2,-2,-2,-3,-2,-2,-2,-2,-3,-2,-2,-3)$.
\end{rem}

\begin{proof}
We use the configuration of Figure~\ref{f:I8} to indicate the embeddings
\begin{figure}[htb]
\begin{center}
\setlength{\unitlength}{1mm}
\includegraphics[height=7cm]{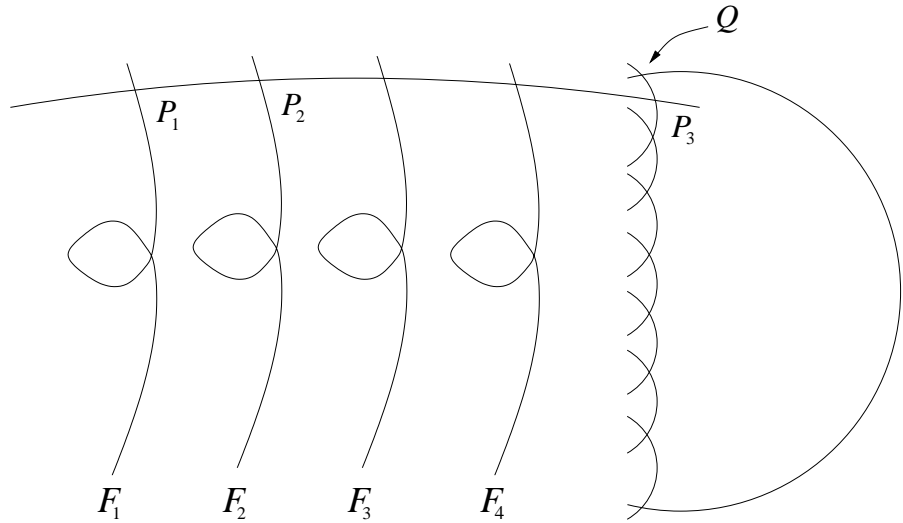}
\end{center}
\caption{\quad A fibration with an $I_8$--fiber}
\label{f:I8}
\end{figure}
given above. First of all, perform the knot surgery in the double node
neighborhood of the fishtail fibers $F_3,F_4$ with isotopic vanishing
cycles and blow up the two double points of the remaining two fishtail
fibers $F_1,F_2$ together with the double point of the
pseudo--section. To get the first embedding, smooth the transverse
intersections $P_2,P_3$ and apply four infinitely close blow--ups at
$Q$, resulting the configuration
\[
(-4,-11,-2,-2,-2,-2,-2,-2,-3,-2,-2,-2). 
\]
One further blow--up of the $(-4)$--sphere provides the first embedding.
If we blow up this sphere as instructed by Figure~\ref{f:inside}, a final
\begin{figure}[htb]
\begin{center}
\setlength{\unitlength}{1mm}
\includegraphics[height=6cm]{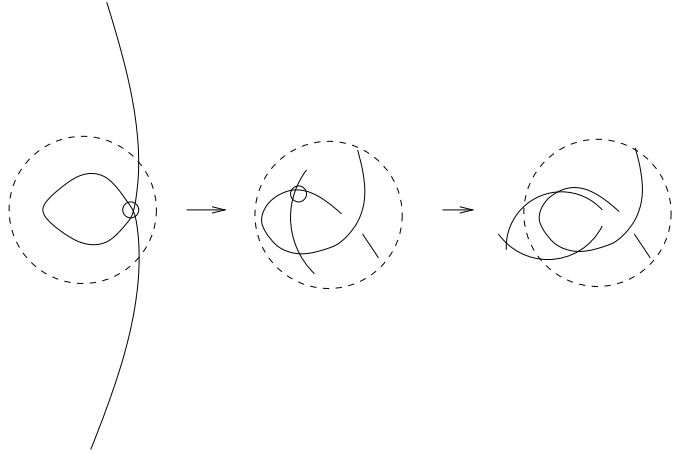}
\end{center}
\caption{\quad Infinitely close blow--ups at the double point of the
fishtail fiber}
\label{f:inside}
\end{figure}
blow--up of the last $(-2)$--sphere in the chain gives the second embedding.

If we smooth the intersections $P_1$ and $P_2$ then eight infinitely
close blow--ups at $Q$, together with a final blow--up on any of the
former fishtail fibers $F_1$ or $F_2$ results the third
embedding. Once again, the last blow--up can be performed as in
Figure~\ref{f:inside}, in which case we need to blow up the other end
of the chain, resulting the fourth embedding. Finally, resolving only
$P_2$, eight infinitely close blow--ups at $Q$, one further blow--up
on the appropriate $(-2)$--sphere in the $I_8$--fiber and one more on
the fishtail passing through $P_1$ gives the fifth configuration.  If
this last blow--up is performed as in Figure~\ref{f:inside}, by
blowing up the last $(-2)$--sphere of the configuration we get the
last promised embedding.  Since in any of the above constructions the
last blow--up provides an exceptional divisor transversally
intersecting the first or last sphere of the configuration, the
complements of the configurations are obviously simply connected.
\end{proof}

Simple Euler characteristic computation
and Freedman's Theorem imply that after rationally blowing down any of the
configurations presented in Proposition~\ref{p:ujconf} we get further 
interesting examples of 4--manifolds homeomorphic to $\cpk \# 5\cpkk$.

\subsubsection{A configuration using the $I_6$--fiber}
\label{ss:i6}

A slightly different procedure can be applied if we start with a
fibration containing an $I_6$--fiber and two pairs of fishtail fibers
with isotopic vanishing cycles. This example was also discovered independently
by Fintushel and Stern \cite{FSuj}. 

\begin{prop}\label{p:i6fibrum}
There is an elliptic fibration $\eegy \to \cp ^1$ with an $I_6$--fiber,
six fishtail fibers $F_1, \ldots , F_6$ and a section.
Furthermore we can assume that the vanishing cycles of $F_1$ and $F_2$
are isotopic, and the vanishing cycles of $F_3$ and $F_4$ are
isotopic.
\end{prop}
\begin{proof}
Start again with the fibration given by the relation
$(a^3b)^3$
and notice that it is equal to 
$ a^6 b ^{a^{-3}}b^2(a^{b^{-1}})^3$.
This expression shows the existence of the required fibration.
\end{proof}

Consider the 4--manifold $V_{K_{n_1},K_{n_2}}$ we get from $\eegy$ by
doing two double node surgeries in the neighborhoods of $F_1,F_2$ and
$F_3,F_4$ respectively, using the knot $K_{n_1}$ for the first and
$K_{n_2}$ for the second surgery. By choosing $K_{n_1},K_{n_2}$ to be
twist knots (as in \cite{FSuj}) we get a pseudo--section $S\subset
V_{K_{n_1},K_{n_2}}$ which is now a sphere with two positive double
points and self-intersection $-1$. Blowing up the two
self--intersections we get a sphere of square $-9$ in
$V_{K_{n_1},K_{n_2}}\# 2 \cpkk$. Using five $(-2)$--spheres of the
$I_6$ fiber, we get a chain of spheres according to the linear
plumbing
\[
(-9,-2,-2,-2,-2,-2),
\]
giving rise to an embedding of $C_{7,1}$ into $V_{K_{n_1},K_{n_2}}\# 2
\cpkk$.  We define our 4--manifolds by rationally blowing down these
copies of $C_{7,1}$.  Simple connectivity of $V_{K_{n_1},K_{n_2}}$
follows from the presence of two different vanishing cycles and the
pseudo--section, while the complement of $C_{7,1}$ in
$V_{K_{n_1},K_{n_2}}\# 2 \cpkk$ is simply connected because there are
two more fishtail fibers in the fibration which we did not use in the
construction. Since $V_{K_{n_1},K_{n_2}}$ is homeomorphic to $\eegy$
and we used two blow--ups to find the above chain of spheres, after
rationally blowing down we get a 4--manifold homeomorphic to $\cpot$.
Recall that $K_n$ denotes the $n$--twist knot (as depicted in
\cite{FSuj}); let $V_n$ denote $V_{K_1, K_n}$. The result of the
rational blow--down of $C_{7,1}\subset V_n \# 2 \cpkk$ will be denoted
by $Q_n$.

\section{Seiberg--Witten computations}
\label{s:harm}
We will prove Theorem~\ref{t:main} by computing Seiberg--Witten
invariants of the 4--manifolds constructed above. We will give
details of the computation for the first construction, resulting the
manifolds $X_n$, very similar ideas work for all the other
manifolds. The argument sketched below is closely modeled on the
argument encountered in \cite{FSuj}. We will finish this section by 
an explicit computation of the Seiberg--Witten invariants of 
the manifolds $Q_n$ constructed in Subsection~\ref{ss:i6}.

It is shown in \cite{FSknot, Sz} that $Y_n$ has two Seiberg--Witten
basic classes $\pm K$, moreover $\vert SW_{Y_n}(\pm K)\vert
=n$. Furthermore, we can choose the sign of $K$ so that it evaluates
on the pseudo--section $S$ as $-1$. Consequently
\[
(K-e_1-\ldots -e_{11}) (u_i)=u_i\cdot u_i +2
\]
for each sphere $u_i$ appearing in the plumbing $C$.  Let $L$ be
the extension of $K\vert _{Y_n-C}$ to $X_n$.  Using the blow--up and
the rational blow--down formula together with the wall--crossing
formula we get 
\begin{prop}\label{p:elso}
The Seiberg--Witten invariant $SW_{X_n}(L)$ is an element  of the set $\{ \pm n
, \pm n \pm 1\}$. Therefore the 4--manifold $X_n$ with $n\geq 2$ admits a
Seiberg--Witten basic class. \qed
\end{prop}
This computation leads us to
\begin{cor}
There exists an exotic smooth 
structure on $\cpot$. 
\end{cor}
\begin{proof}
Since the Seiberg--Witten function is a diffeomorphism invariant for
manifolds with $b_2^+=1$ and $b_2^-\leq 9$, and by the existence of a
positive scalar curvature metric we have $SW_{\cpot}\equiv 0$, we get
that $X_n$ is not diffeomorphic to $\cpot$, hence the corollary
follows.
\end{proof}

Since $Y_n$ has exactly two basic classes, the same computation as
above actually shows
\begin{lem}\label{l:masodik}
The Seiberg--Witten function of $X_n$ takes its values in a subset of
$\{ 0, \pm 1 ,\pm n , \pm n \pm 1\}$, and for $n\geq 3$ there are
exactly two basic classes $\pm L$ with Seiberg--Witten values in $\{
\pm n , \pm n \pm 1\}$.  \qed
\end{lem}

\begin{proof}[Proof of Theorem~\ref{t:main}]
Combining Proposition~\ref{p:elso} with Lemma~\ref{l:masodik} it
follows that $X_n$ and $X_{n+3k}$ are not diffeomorphic once $n\geq 2$
and $k> 0$. This observation proves the existence of infinitely many
distinct smooth structures on $\cpot$.  The blow--up formula and the
fact that for $n\geq 3$ there are only two basic classes of $X_n$ with
Seiberg--Witten values in the set $\{ \pm n , \pm n \pm 1\}$ show
that for $n\geq 3$ the manifold $X_n$ is actually minimal.
\end{proof}

The argument above was sufficient for proving Theorem~\ref{t:main},
but with some additional work the complete Seiberg--Witten invariants
of the 4--manifolds encountered above can be determined. We 
demonstrate this for the 4--manifolds $Q_n$ defined in
Subsection~\ref{ss:i6} and prove

\begin{thm}\label{t:qn}
For $n \geq 1$ the 4--manifold $Q_n$ admits exactly two basic classes
$\pm L$ and $SW_{Q_n}(\pm L)=\pm n$. Consequently the manifolds $Q_n$ are
all minimal and pairwise nondiffeomorphic.
\end{thm}

The heart of the argument is to find a simple way to relate the
Seiberg--Witten invariants of $Q_n$ to those of $V_n$.
As a stepping stone we will need the following construction.

Start with the fibration $\eegy \to S^2$ provided by
Proposition~\ref{p:i6fibrum}. Instead of doing the double node
surgery, blow up the 4--manifold $\eegy$ twice and in the two new
$\cpkk$'s choose embedded spheres representing twice the generator of
$H_2(\cpkk ; \bfz )$.  By tubing these two $(-4)$--spheres to a fixed
section of $\eegy \to S^2$ we get a $(-9)$--sphere, which, together
with five $(-2)$--spheres of the $I_6$--fiber gives rise to an
embedded copy of $C_{7,1}$ in $( \eegy ) \# 2 \cpkk=\cpk \# 11 \cpkk$.
Let $R$ denote the 4--manifold we get by rationally blowing down this
copy of $C_{7,1}$.

\begin{prop}
The  Seiberg--Witten invariant $SW_R$ is identically zero.
\end{prop}
\begin{proof}
Note that $b_2^-(R)=5$, hence the Seiberg--Witten function $SW_R$ is
well--defined.  Let $D\subset \eegy$ denote the tubular neighborhood
of the chosen section and the chain of five $(-2)$--spheres in the
$I_6$--fiber. Notice that $\partial D=S^3$. By performing the
blow--ups and the rational blow--down process in $D$ (resulting in a
negative definite 4--manifold $W$), we get a decomposition of $R$ as
$P\# W$. Since $P\# 6 \cpkk=\eegy$, the blow--up formula and
$SW_{\eegy }\equiv 0$ imply that $SW_P\equiv 0$. Now the usual gluing
formula along $S^3$ implies the result.
\end{proof}

Notice that by the construction of $V_n$ there is a natural bijection
\[
{\mathcal {F}}\colon H_2 (\eegy ; \bfz ) \to H_2 (V_n ; \bfz )
\] 
mapping the chosen section of $\eegy \to S^2$ to the pseudo--section
in $V_n$. The map ${\mathcal {F}}$ induces 
a natural extension  to the double
blow--ups 
\[
{\mathcal {F}}'\colon H_2 (\cpk \# 11 \cpkk ; \bfz ) \to H_2 (V_n \#
2 \cpkk ; \bfz ).
\] 
In these double blow--ups we have found copies of $C_{7,1}$; 
%denote by $C'\subset \cpk \# 11\cpkk$ and $C''\subset V_n \# 2\cpkk$.
it follows from the constructions of these submanifolds that
${\mathcal {F}}'$ maps the homology classes of the
chains of spheres into each other.

In addition, homology classes of $R$ (resp. $Q_n$) can be naturally
constructed from homology classes of $\cpk \# 11 \cpkk $ (resp. $V_n\#
2\cpkk$) by  appropriately extending them to the
rational blow--down.  In particular, the map ${\mathcal {F}}'$ gives
rise to a bijection
\[
{\mathcal {F}}_1\colon H_2 (R ; \bfz ) \to H_2 (Q_n ; \bfz ).
\]

Let $K\in H^2 (V_n ; \bfz )$ be a characteristic element.  For odd
integers $a,b$ we get extensions $K+aE_1 +bE_2\in H^2 (V_n \# 2\cpkk ;
\bfz )$, where $E_i$ denote the Poincar\'e duals of the exceptional
divisors of the blow--ups. Suppose that the restriction of $K+aE_1
+bE_2$ to $V_n \# 2\cpkk - C_{7,1}$ extends to a characteristic
cohomology class to $Q_n$ and denote this extension by $K(a,b)$.
Suppose furthermore that the formal dimension $d(K+aE_1+bE_2)$ of the
Seiberg--Witten moduli space on $V_n \# 2\cpkk$ corresponding to
$K+aE_1+bE_2$ is nonnegative. 
\begin{lem} \label{l:ut}
Let $K, a, b$ be chosen as above. Then
\[
SW_{Q_n}(K(a,b))-SW_{V_n}(K)= SW_{R}(f_1(K(a,b)))-
SW_{\eegy}(f(K)),
\]
where $f$ and $f_1$ are duals of ${\mathcal {F}}$ and ${\mathcal
{F}}_1$.
\end{lem}
\begin{proof}
Since the blow--up, wall--crossing and rational blow--down formulae
involve only homological computations, and ${\mathcal {F}}'$
identifies the two copies of $C_{7,1}$, the lemma follows.
\end{proof}

\begin{proof}[Proof of Theorem~\ref{t:qn}]
Let $L\in H^2 (Q_n ; \bfz )$ be a characteristic element with $SW_{Q_n
}(L)\neq 0$. By the rational blow--down formula there is a class
$K+aE_1+bE_2\in H^2 (V_n \# 2\cpkk ; \bfz )$ with 
\[
SW_{Q_n}(L)=SW_{V_n \# 2\cpkk }(K+aE_1+bE_2)
\]
where the right--hand side is taken in the appropriate chamber.  In
particular, $L=K(a,b)$ for some $K\in H^2 (V_n ; \bfz )$ and
$d(K+aE_1+bE_2)\geq 0$.  Since $SW_{R}\equiv 0$ and $SW _{\eegy}
\equiv 0$, Lemma~\ref{l:ut} implies that $SW_{V_n}(K)\neq 0$.  On the
other hand, the Seiberg--Witten invariants of $V_n$ are known
\cite{FSknot}, hence it follows that $K=\pm T, \pm 3T$ where $T$ is
the Poincar\'e dual of the fiber. Since $d(T)=d(3T)=0$, it follows
that $a=\pm 1$ and $b=\pm 1$.  A simple homological computation shows
that in the family $\{ \pm T \pm E_1\pm E_2, \ \pm 3T \pm E_1\pm E_2\}
\subset H^2 (V_n \# 2\cpkk ; \bfz )$ there are only two cohomology
classes --- which are equal to $\pm (3T-E_1-E_2)$ --- admitting
extensions to $Q_n$. Since $SW_{V_n}(\pm 3T)=\pm n$ the theorem
follows from Lemma~\ref{l:ut}.
\end{proof}

\end{document}